\documentclass{daj}

\usepackage{amsmath, amssymb, amsfonts, amsthm, bbm, mathtools, mathrsfs}
\usepackage{verbatim,enumitem,graphics}

\newtheorem{theorem}{Theorem}[section]
\newtheorem{question}[theorem]{Question}
\newtheorem{conjecture}[theorem]{Conjecture}

\theoremstyle{definition}
\newtheorem{definition}[theorem]{Definition}
\newtheorem{example}[theorem]{Example}
\newtheorem{remark}[theorem]{Remark}

\newcommand{\cF}{\mathcal{F}}

\newcommand{\vx}{\Vec{x}}
\newcommand{\vy}{\Vec{y}}

\def\l{\ell }
\def\e{\varepsilon}

\dajAUTHORdetails{%
  title = {Color Avoidance for Monotone Paths}, 
  author = {Eion Mulrenin, Cosmin Pohoata, and Dmitrii Zakharov},
  plaintextauthor = {Eion Mulrenin, Cosmin Pohoata, Dmitrii Zakharov},
  keywords = {Ramsey theory, monotone paths in hypergraphs},
}   

\dajEDITORdetails{%
   year={2025},
   number={23},
   received={8 January 2025},   
   published={7 October 2025},  
   doi={10.19086/da.145175},       
}   

\begin{document}

\begin{frontmatter}[classification=text]

\title{Color Avoidance for Monotone Paths} 

\author[eion]{Eion Mulrenin}
\author[cosmin]{Cosmin Pohoata\thanks{Supported by NSF Award DMS-2246659.}}
\author[dima]{Dmitrii Zakharov\thanks{Supported by the Jane Street Graduate Fellowship.}}

\begin{abstract}
    In 2014, Moshkovitz and Shapira determined the tower height for hypergraph Ramsey numbers of tight monotone paths. We address the color-avoiding version of this problem in which one no longer necessarily seeks a monochromatic subgraph, but rather one which {\it avoids} some colors. This problem was previously studied in uniformity two by Loh and by Gowers and Long. 
    
    We show, in general, that the tower height for such Ramsey numbers requires one less exponential than in the usual setting.
    The transition occurs at uniformity three, where the usual Ramsey numbers of monotone paths of length $n$ are exponential in $n$, but the color-avoiding Ramsey numbers turn out to be polynomial.
\end{abstract}
\end{frontmatter}

\section{Introduction}

\subsection{Background}
It was noticed a little more than a decade ago \cite{EM13, FPSS12} that two seminal theorems of Erd\H{o}s and Szekeres from their 1935 paper \cite{ES35}---the \textit{monotone subsequence theorem} and the \textit{cups-caps lemma}---can be cast as Ramsey theoretic statements about monotone paths in two- and three-uniform hypergraphs, respectively.
A \textit{tight monotone path of length $n$} in a $k$-uniform hypergraph on $[N]$ is the set of $n$ edges which contain $k$ consecutive vertices from $v_1 < v_2 < \dots < v_{n+k-1}$ in $[N]$, and we will denote by $MS_k(n; r)$ the least integer $N$ such that every $r$-coloring of $[N]^{(k)}$ contains a tight \textit{monochromatic} monotone path of length $n$, i.e. a monotone path in which every edge receives the same color.
(Here and throughout, for a natural number $N$, we will denote by $[N]$ the set $\{1,2,\dots,N\}$, and for a set $X$ and a natural number $k$, by $X^{(k)}$ the set $\{e \subseteq X: |e| = k\}$.)
In this language, the theorems of Erd\H{o}s and Szekeres read as follows.

\begin{theorem}[Erd\H{o}s--Szekeres, 1935 \cite{ES35}]
    For all positive integers $n \geq 3$ and $r \geq 2$, we have
    \begin{equation*}
        MS_2(n; 2) = n^2 + 1.
    \end{equation*}
    Also, for $n \geq 4$,
    \begin{equation*}
        MS_3(n; 2) = \binom{2n}{n} + 1.
    \end{equation*}
\end{theorem}

In \cite{FPSS12}, the authors identified this unifying framework and asked for sharp estimates on $MS_k(n; r)$ for larger values of $k$ and $r$.
Soon after, the problem was almost entirely resolved by Moshkovitz and Shapira \cite{MS14}, who established a connection between Ramsey numbers of monotone paths and certain types of higher-order down-sets as well as high-dimensional integer partitions.
They proved an exact equality between $MS_k(n; r)$ and the number of these down sets which yielded the following quantitative estimates.
To state them, we recursively define the \textit{tower function} $T_k(m)$ by $T_0(m) = m$ and $T_{k+1}(m) = 2^{T_k(m)}$; note that the subscript denotes the number of 2's in the tower.

\begin{theorem}[Moshkovitz--Shapira, 2014 \cite{MS14}]\label{MS}
    For all $k \geq 3$, $r \geq 2$, and $n \geq 2$,
    \begin{equation*}
        MS_k(n; r) \leq T_{k-2}(2 n^{r-1}).
    \end{equation*}
    Moreover, for $n$ sufficiently large,
    \begin{equation*}
        MS_k(n; r) \geq T_{k-2}(n^{r-1} / 2 \sqrt{r}).
    \end{equation*}
\end{theorem}

\begin{remark}
    The second and third authors have recently shown~\cite{PZ24} that this lower bound gives the correct order of magnitude, i.e., $MS_k(n;r) = T_{k-2}(\Theta(n^{r-1}/\sqrt{r}))$ when $r$ is large. 
    See also \cite{falgas2023dedekind}.
\end{remark}

On the other hand, another line of research which has received considerable attention concerns Ramsey problems in which the desired structure now \textit{avoids} some number of colors instead of using only one color.
Note that this generalizes the classical problem since a subgraph in a $r$-colored graph which avoids $r-1$ colors is necessarily monochromatic.
This problem was introduced in a 1965 paper by Erd\H{o}s, Hajnal, and Rado~\cite{EHR65}, and was improved several years later by Erd\H{o}s and Szemer\'edi \cite{ES72}.
For cliques, the problem has seen a number of new results in the past several years \cite{ACMMM24, CFHMSV24, DGHY23}.

In the context of monotone paths, such color-avoiding Ramsey problems were first considered by Loh \cite{Loh15} in 2015. 
More precisely, he studied the parameter $A_2(n; 3,2)$ which we shall define in greater generality as follows.

\begin{definition}
    For positive integers, $n \geq k$ and $r > s$, let $A_k(n; r, s)$ denote the least $N$ such that every $r$-coloring of $[N]^{(k)}$ contains a tight monotone path on $n$ edges which uses at most $s$ colors.
\end{definition}

A trivial bound for $A_{2}(n;3,2)$ comes from the observation that a $3$-coloring of $[N]^{(2)}$ can be regarded as a $2$-coloring of $[N]^{(2)}$ if one pretends that two of the colors are the same, and that a monochromatic path in the latter setting is a path which avoids the usage of a color in the original $3$-coloring. The Erd\H{o}s--Szekeres theorem then implies that
$$A_2(n; 3, 2) \leq MS_{2}(n;2) = n^2+1.$$
In \cite{Loh15}, Loh noted that $A_{2}(n;3,2) \geq n^{3/2}$ and used the proof of the Erd\H{o}s--Szekeres theorem due to Seidenberg \cite{Seid59} together with the triangle removal lemma of Ruzsa and Szemer\'edi \cite{RS78} to show that $A_{2}(n;3,2)=o(n^2)$ holds. In a remarkable paper \cite{GL21}, Gowers and Long further showed that there exists an absolute constant $\e > 0$ such that $A_2(n; 3, 2) \leq n^{2-\e}$ and conjectured that $A_{2}(n;3,2) = \Theta(n^{3/2})$.

\subsection{Our results}

A more general first estimate on $A_{k}(n;r,s)$ can be obtained in a similar fashion from Theorem \ref{MS} by noticing that an $r$-coloring of $[N]^{(k)}$ can be regarded as a $\lceil r/s \rceil$-coloring of $[N]^{(k)}$ if one pretends that $\lceil r/s \rceil$ disjoint groups of $s$ colors each represent a single color, and that a monochromatic tight path in the latter setting is a tight path which uses at most $s$ colors in the original $r$-coloring. For example, when $k=3$, $r=3$, and $s=2$, this gives
$$A_3(n; 3, 2) \leq MS_3(n; 2) = \binom{2n}{n} + 1.$$
In general, for every $n \geq k \geq 3$, $r > s \geq 1$,
$$A_k(n; r, s) \leq MS_k(n; \lceil r/s \rceil) \leq T_{k-2}(2n^{\lceil r/s \rceil}).$$

Our main result is an improvement on the height of the tower in these upper bounds by one exponential, in the case when $k \geq 3$ and $s > r/2$. For example, when $k=3$, we show the following:

\begin{theorem}\label{main3}
    There exists an absolute constant $C > 0$ such that for all positive integers $n \geq 3$ and $r$, $s$ with $s > r/2$, 
    \begin{equation*}
        A_3(n; r, s) \leq n^{C \binom{r}{s}^2 \log \binom{r}{s}}.
    \end{equation*}
\end{theorem}

Moreover, for the particular combination $r=3$ and $s=2$, our proof gives the stronger bound $A_{3}(n;3,2) \leq n^{9}$. Qualitatively speaking, Theorem \ref{main3} shows that both $A_2(n;r,s)$ and $A_3(n;r,s)$ are polynomial in $n$ for $s > r/2$. We also show that for uniformity $k \geq 4$, the tower height continues to require one less exponential than the trivial bound from the theorem of Moshkovitz and Shapira.

\begin{theorem}\label{maink}
    There exists an absolute constant $C > 0$ such that for all positive integers $n \geq k \geq 3$ and $r$, $s$ with $s > r/2$, we have
    \begin{equation*}
        A_k(n; r,s) \leq T_{k-3}\left(n^{C \binom{r}{s}^2 \log \binom{r}{s}}\right).
    \end{equation*}
\end{theorem}

\subsection{Organization}
The remainder of the paper is organized as follows.
In Section~\ref{prelim}, we give an overview of some preliminaries which will be used extensively in our proofs.
In Section~\ref{s3}, we focus on the problem in uniformity three, i.e., on $A_3(n;r,s)$ and Theorem~\ref{main3}.
We give first a proof that $A_3(n;3,2) \leq n^9$ followed by the general proof, as it is a direct generalization of the argument for $r=3$ and $s=2$. 
In Section~\ref{sk}, we prove Theorem~\ref{maink} by adapting the method of Section~\ref{s3}, and in Section~\ref{scharacterization} we give a characterization of $A_k(n;r,s)$ in terms of what we call ``generalized $S$-increasing sequences". 
Finally, in Section~\ref{conclusion}, we make some closing remarks and state some open problems.

\section{Preliminaries}\label{prelim}

In this subsection, we would like to first revisit the fact that
\begin{equation} \label{cupcap}
        MS_3(n; 2) \leq \binom{2n}{n} + 1
    \end{equation}
    holds for every $n \geq 4$. In \cite{ES35}, Erd\H{o}s and Szekeres originally proved this result in the language of cups and caps by establishing a recursive inequality in the spirit of the classical recursive inequality for diagonal Ramsey numbers (also originally from \cite{ES35}). In \cite{CK71} and \cite{MS14}, however, Chv\'atal--Koml\'os and Moshkovitz--Shapira note that it is also possible to prove \eqref{cupcap} by appropriately extending the Seidenberg argument \cite{Seid59}. Our proof of Theorem \ref{main3} (and that of Theorem \ref{maink}) will also adopt this perspective, so for the sake of self-containment we will now recall this beautiful idea below. 

   \subsection{Seidenberg's argument in uniformity three.} \label{Seidenberg} Let $\Delta: [N]^{(3)} \to \{\operatorname{red},\operatorname{blue}\}$ be a two-coloring of the triples of $[N]$ with no tight monotone monochromatic path of length $n$.
    For each $u < v$, define the vector
    \begin{equation*}
        P(u,v) = \begin{bmatrix}
            \l_{\text{red}}(u,v) + 1\\
            \l_{\text{blue}}(u,v) + 1
        \end{bmatrix},
    \end{equation*}
    where $\l_{\text{red}}(u,v)$ and $\l_{\text{blue}}(u,v)$ denote the lengths of the longest tight monotone paths in red and blue, respectively, which end at $u, v$.
    Since there is no tight monotone path of length $n$, note that the point $P(u,v)$ lies in the grid $[n]^2$. Now, consider the map
    \begin{equation*}
        \Theta: v \longmapsto D(v) = \{\vx \in [n]^2: \vx \leq P(u,v) \text{ for some } u < v\},
    \end{equation*}
    which assigns elements of $[N]$ down-sets of $[n]^2$. Here the $\leq$ relation on $[n]^2$ denotes the usual ``less than or equal to in both coordinates" partial order on $[n]^2$. One can check that $D(v)$ is indeed a down-set: for every $\vx \in D(v)$ we have that $\vx' \in D(v)$ for every $\vx' \leq \vx$. 
    
    The main observation is that the map $\Theta$ is injective: if $u < v$ and $D(u) = D(v)$, we have $P(u,v) \in D(v) = D(u)$, so by definition $P(u,v) \leq P(w,u)$ for some $w < u$; but this is a contradiction since the longest path ending at $u, v$ in color $\Delta(w,u,v)$ is strictly longer than the one ending at $w,u$.
    Since the number of down-sets in $[n]^2$ equals the number of antichains in $[n]^2$ which is $\binom{2n}{n}$, we have that $N \leq \binom{2n}{n}$. We refer to \cite{MS14} for more details regarding the latter claims. 


    The main idea behind Theorem \ref{main3} is that in order to upper bound $A_{3}(n;r,s)$ one can construct an injection of $[N]$ into a significantly smaller space than the collection of down-sets of $[n]^2$. To define this map, we will require the concept of a dominating set in a majority tournament. 
    
\subsection{Majority tournaments}

A {\it{tournament}} is a directed graph $T=(V,E)$ where every two vertices $u,v \in V$ are connected by one edge, i.e. $|E \cap \left\{(u,v),(v,u)\right\}|=1$. If the edge is $(u,v)$ (denoted by $u \to v$) we say that $u$ {\it{dominates}} $v$. 

\begin{definition}
    For an integer $K \geq 2$, a tournament $T$ is called a $K$-\textit{majority tournament} if there exist $2K-1$ linear (total) orders $\prec_1, ..., \prec_{2K-1}$ on the set of vertices of $T$ such that $u \to v$ in $T$ if and only if $v \prec_i u$ holds for at least $K$ values of $i \in \left\{1,\ldots,2K-1\right\}$. 
\end{definition}

These types of tournaments have been studied in social science~\cite{Mcg53} as well as combinatorics, where previous researchers studied extremal problems regarding the sizes of induced transitive subtournaments \cite{MSW11} and dominating sets \cite{ABKKW06}. If $T$ is a tournament and $X,Y$ are subsets of vertices of $T$, we say that $X$ {\it{dominates}} $Y$ if for every vertex $y$ in $Y \setminus X$ there exists an $x \in X$ such that $x$ dominates $y$. If $Y$ is the full set of vertices of $T$, then we call $X$ a {\it{dominating set}} of the tournament $T$. The {\it{domination number}} $f(T)$ of $T$ is the smallest cardinality of a dominating set of $T$. It is a classical fact (often attributed to Erd\H{o}s, see for example \cite{Moon68}) that every tournament $T$ on $n$ vertices satisfies $f(T) \leq \lceil \log_{2}{n} \rceil$. It is also known \cite{AS16} that random tournaments on $n$ vertices do not contain any smaller dominating sets, so this estimate cannot be improved upon in general. Nevertheless, it turns out that $K$-\textit{majority tournaments} behave quite differently than random tournaments when it comes to their domination number $f(T)$. In \cite{ABKKW06}, Alon, Brightwell, Kierstead, Kostochka, Winkler proved the following result:

\begin{theorem} \label{kmaj}
If $T$ is a $2$-majority tournament, then 
$$f(T) \leq 3.$$
In general, for every $K \geq 2$, if $T$ is a $K$-majority tournament, then
$$f(T) \leq C K \log K$$
holds for an absolute constant $C>0$. 
\end{theorem}

It turns out it is this phenomenon which is primarily responsible for the drop in the tower height in our estimates from Theorem \ref{main3} and Theorem \ref{maink} when compared to the estimate from Theorem \ref{MS}. 


\section{Proof of Theorem~\ref{main3}}\label{s3}

\subsection{The proof for $r=3$, $s=2$.}

We first prove Theorem~\ref{main3} for $r=3$ and $s=2$, in order to showcase the main idea in the clearest way possible. Our proof for this regime actually gives the stronger bound $A_3(n;3,2) \leq n^9$.

    Let $\Delta: [N]^{(3)} \to [3]$ be a 3-coloring of the triples of $[N]$ in which every tight monotone path of length $n$ receives all three colors.
    Define for each $u < v$ in $[N]$ the vector
    \begin{equation*}
        P(u,v) = \begin{bmatrix}
            \l_{1,2}(u,v) + 1\\
            \l_{2,3}(u,v) + 1\\
            \l_{1,3}(u,v) + 1
        \end{bmatrix} \in [n]^3,
    \end{equation*}
    where $\l_{i,j}(u,v)$ is the length of the longest tight monotone path ending at $u, v$ which uses only colors $i$ and $j$.
    
    Like in the argument from Subsection \ref{Seidenberg}, let us consider the map
    \begin{eqnarray*}
        \Theta_1: v \longmapsto S(v) = \{P(u,v) \in [n]^3: u < v\},
    \end{eqnarray*}
    associating elements of $[N]$ with subsets of $[n]^3$. 
    To define an injection from $[N]$ into a space of size smaller than ${2n \choose n}$, we define a 2-majority tournament on the elements of $S(v)$ with linear orders $\prec_1, \prec_2, \prec_3$ as follows.
    If $\vx = (x_1,x_2,x_3), \vy = (y_1,y_2,y_3) \in S(v)$, set
    \begin{align*}
        & 1. \hspace{0.25cm} \Vec{x} \prec_1 \Vec{y} \text{ if } \begin{cases}
            x_1 < y_1; \text{ or}\\
            x_1=y_1, x_2 < y_2; \text{ or}\\
            x_1=y_1, x_2=y_2, x_3 < y_3
        \end{cases}\\
        & 2. \hspace{0.25cm} \Vec{x} \prec_2 \Vec{y} \text{ if } \begin{cases}
            x_2 < y_2; \text{ or}\\
            x_2=y_2, x_3 < y_3; \text{ or}\\
            x_2=y_2, x_3=y_3, x_1 < y_1
        \end{cases}\\
        & 3. \hspace{0.25cm} \Vec{x} \prec_3 \Vec{y} \text{ if } \begin{cases}
            x_3 < y_3; \text{ or}\\
            x_3=y_3, x_1 < y_1; \text{ or}\\
            x_3=y_3, x_1=y_1, x_2 < y_2
        \end{cases}
    \end{align*}
    
    In words, $\prec_i$ is the lexicographic order on $S(v) \subseteq [n]^3$ in which the $i$th coordinate is the most important, followed by the $(i+1)$-st and then the $(i+2)$-nd (with addition modulo 3).
    By Theorem \ref{kmaj}, $S(v)$ has a total dominating set $D(v)$ of size at most three, i.e., for every $\vy \in S(v)$, there exists $\vx \in D(v)$ with $\vx \to \vy$ in the tournament defined above. 
    
    Now consider a second map
    \begin{equation*}
        \Theta_2: v \longmapsto D(v) \subseteq S(v).
    \end{equation*}
    Since $D(v) \subseteq S(v) \subseteq [n]^3$, note that $\Theta_{2}$ is a map from $[N]$ to the set of subsets of size at most $3$ of $[n]^3$. We claim that $\Theta_2$ is also an injection.
    Suppose, on the contrary, that $u < v$ but $D(u) = D(v)$.
    Then each element $\Vec{x} \in D(v) = D(u) \subseteq S(u)$ by definition has the form $\Vec{x} = P(w,u)$ for some $w < u$;
    but $\Vec{x} = P(w,u) \prec_i P(u,v)$ for at least the two coordinates $i$ corresponding to paths in which the color $\Delta(w,u,v)$ may be used, and so $P(u,v) \to \vx$ for all $\vx \in D(v)$, contradicting the fact that $D(v)$ totally dominates $S(v)$ in the 2-majority tournament defined by the $\prec_i$'s.

    Therefore, 
    we have $N \leq \binom{n^3}{3}+\binom{n^3}{2} +\binom{n^3}{1}+\binom{n^3}{0} \le n^9$ as promised. \qed


\subsection{General $r$ and $s$.}

Now we prove Theorem~\ref{main3} for all combinations $r, s$ with $s > r/2$.
The proof is a straightforward generalization of the one given in the previous subsection.

    Let $\Delta: [N]^{(3)} \to [r]$ be an $r$-coloring of the triples of $[N]$ in which every tight monotone path on $n$ edges receives at least $s+1$ colors (i.e., no such path uses at most $s$ colors).
    As before, define for each $u < v$ the vector
    \begin{equation*}
        P(u,v) = \langle 1 + \l_S(u,v)\rangle_{S \in [r]^{(s)}} \in [n]^{\binom{r}{s}},
    \end{equation*}
    where $\l_S(u,v)$ denotes the length of the longest monotone path ending at $u$ and $v$ which only uses colors from the set $S \in [r]^{(s)}$,
    and consider the map
    \begin{equation*}
        \Theta_1: v \longmapsto S(v) = \{P(u,v) \in [n]^{\binom{r}{s}}: u < v\}
    \end{equation*}
    associating elements of $[N]$ with subsets of $[n]^{\binom{r}{s}}$.

    Let $K$ be the positive integer satisfying $\binom{r}{s} = 2K-1$ if $\binom{r}{s}$ is odd and $\binom{r}{s} = 2(K-1)$ if $\binom{r}{s}$ is even.
    For each $v \in [N]$, we define a $K$-majority tournament on the elements of $S(v)$ with linear orders $\prec_1, \prec_2, \dots, \prec_{\binom{r}{s}}$ as follows.
    If $\vx = (x_1,x_2, \dots, x_{\binom{r}{s}}), \vy = (y_1,y_2, \dots, y_{\binom{r}{s}}) \in S(v)$ and $1 \leq i \leq \binom{r}{s}$ (with addition in the indices considered modulo $\binom{r}{s}$), set  
    \begin{equation*}
        \vx \prec_i \vy \text{ if } \begin{cases}
            x_i < y_i; \text{ or}\\
            x_i=y_i, x_{i+1} < y_{i+1}; \text{ or}\\
            x_i=y_i, x_{i+1} = y_{i+1}, x_{i+2} < y_{i+2}; \text{ or}\\
            \vdots\\
            x_i=y_i, x_{i+1} = y_{i+1}, x_{i+2} = y_{i+2}, \dots, x_{i+\binom{r}{s}-2} = y_{i+\binom{r}{s}-2}, x_{i+\binom{r}{s}-1} < y_{i+\binom{r}{s}-1}
        \end{cases}
    \end{equation*}
    Note that if $\binom{r}{s} = 2K-1$ is odd then this defines a $K$-majority tournament on $S(v)$.
    If $\binom{r}{s} = 2(K-1)$ is even, then we will force the above construction to be a $K$-majority tournament by adding an additional arbitrary linear order $\prec_{*}$ with no relation to the others.
    By Theorem \ref{kmaj}, each $S(v)$ has a total dominating set $D(v)$ in the tournament defined above, i.e., for every $\vy \in S(v)$, there exists $\vx \in D(v)$ with $\vx \to \vy$, with $|D(v)| \le C K \log K$.

    Now consider a second map
    \begin{equation*}
        \Theta_2: v \longmapsto D(v) \subseteq S(v).
    \end{equation*}
    We argue that $\Theta_2$ is again an injection.
    To see this, assume otherwise that for some $u < v$ in $[N]$, $D(u) = D(v)$.
    Each $\vx \in D(v) = D(u) \subseteq S(u)$ has the form $P(w,u)$ for some $w < u$.
    Now $\vx = P(w,u) \prec_i P(u,v)$ in all of the $\binom{r-1}{s-1}$ coordinates $i \in [\binom{r}{s}]$ which correspond to subsets $S \subseteq [r]$ containing the color $\Delta(w,u,v)$.
    By our hypothesis that $s > r/2$, we have
    \begin{equation*}
        2(K-1) \leq \binom{r}{s} = \frac{r}{s} \cdot \binom{r-1}{s-1} < 2 \cdot \binom{r-1}{s-1},
    \end{equation*}
    which gives (regardless of parity, since all values are integers) that 
    \begin{equation*}
        \binom{r-1}{s-1} \geq K.
    \end{equation*}
    Thus, each such $\vx$ has $\vx \prec_i P(u,v)$ for at least $K$ of the orders $\prec_i$, $1 \leq i \leq \binom{r}{s}$, and so $P(u,v) \to \vx$ in the $K$-majority tournament on $S(v)$ for all $\vx \in D(v)$, contradicting the fact that $D(v)$ is a total dominating set of $S(v)$.

    Therefore, since each $D(v)$ has size at most $C K \log K = C' \binom{r}{s} \log \binom{r}{s}$, we have
    \begin{equation*}
        N \leq \binom{n^{\binom{r}{s}}}{\le C' \binom{r}{s} \log \binom{r}{s}} \leq n^{C' \binom{r}{s}^2 \log \binom{r}{s}}.
    \end{equation*}
    \qed

\vspace{-10mm}

\section{Proof of Theorem~\ref{maink}}\label{sk}

The argument is an extension of the one in uniformity three, but we no longer define an injection from the singletons into the subsets of $[n]^{\binom{r}{s}}$ of size at most $C \binom{r-1}{s-1} \log \binom{r-1}{s-1}$.
Instead, we define such a map on $[N]^{(k-2)}$ and show that it contains no 2-edge monotone path in which both edges are mapped to the same set, or, in other words, we define a coloring of $[N]^{(k-2)}$ with $n^{C \binom{r}{s} \binom{r-1}{s-1} \log \binom{r-1}{s-1}}$ colors and no monochromatic monotone path of length two.
This allows us to use Theorem~\ref{MS} to get an improved upper bound for $A_k(n;r,s)$ over the trivial $MS_k(n;\lceil r/s \rceil)$.
Similar ideas were used in \cite{FPSS12, MS14} to upper bound $MS_k(n,r)$.\\

\noindent \textit{Proof of Theorem~\ref{maink}}.
    Let $\Delta: [N]^{(k)} \to [r]$ be an $r$-coloring of the $k$-subsets of $[N]$ in which every tight monotone path on $n$ edges receives at least $s+1$ colors.
    For $v_1 < v_2 < \dots < v_{k-1}$, define
    \begin{equation*}
        P(v_1,v_2,...,v_{k-1}) =
            \langle 1 + \l_S(v_1, v_2, ...,v_{k-1})\rangle_{S \in [n]^{(r)}} \in [n]^{\binom{r}{s}},
    \end{equation*}
    where as above $\l_S(v_1,...,v_{k-1})$ denotes the length of the longest tight monotone path ending at $v_1, ..., v_{k-1}$ which only uses colors in the set $S \in [r]^{(s)}$.
    
    Consider the function
    \begin{align*}
         \Theta_1: \{v_1,...,v_{k-2}\} \longmapsto S(v_1,...,v_{k-2}) = \{P(u,v_1,...,v_{k-2}) \in [n]^{\binom{r}{s}}: u < v_1\}
    \end{align*}
    mapping elements of $[N]^{(k-2)}$ into subsets of $[n]^{\binom{r}{s}}$.
    We can define the exact same $K$-majority tournament as above on each $S(v_1,...,v_{k-2})$ and obtain a total dominating set $D(v_1,...,v_{k-2}) \subseteq S(v_1,...,v_{k-2})$ of size at most $C K \log K$.
    With this, again consider a second map
    \begin{equation*}
        \Theta_2: \{v_1,...,v_{k-2}\} \longmapsto D(v_1,...,v_{k-2}) \subseteq S(v_1,...,v_{k-2}).
    \end{equation*}
    Now let $S \subseteq [n]^{\binom{r}{s}}$ have size $|S| \le C K \log K$.
    We claim that $\Theta_2^{-1}(S)$ has no tight $(k-2)$-uniform monotone path of length two.
    Indeed, if it had such a path on vertices $v_1 < v_2 < ... < v_{k-1}$, then each $\Vec{x} \in D(v_2,...,v_{k-1}) = S = D(v_1,...,v_{k-2}) \subseteq S(v_1,...,v_{k-2})$ by definition has $\Vec{x} = P(v_0,v_1,...,v_{k-2})$ for some $v_0 < v_1$; but then
    \begin{equation*}
        \Vec{x} = P(v_0,...,v_{k-2}) \prec_i P(v_1,...,v_{k-1}) \in S(v_2,...,v_{k-1})
    \end{equation*}
    in the $\binom{r-1}{s-1} \geq K$ coordinates $i$ corresponding to paths which are allowed to contain the color $\Delta(v_0,...,v_{k-1})$,
    and so $P(v_1, \dots, v_{k-1}) \to \vx$ for all $\vx \in D(v_2, \dots, v_{k-1})$, contradicting the fact that $D(v_2,...,v_{k-1})$ was chosen as a total dominating set for $S(v_2,...,v_{k-1})$.
    
    In other words, then, the partition
    \begin{equation*}
        [N]^{(k-2)} = \bigcup_{S \subseteq [n]^{\binom{r}{s}}, |S| \le C K \log K} \Theta_2^{-1}(S)
    \end{equation*}
    defines an $\binom{n^{\binom{r}{s}}}{\le C K \log K}$-coloring of the $(k-2)$-tuples of $[N]$ with no monochromatic monotone path of length two.
    Thus, by Theorem \ref{MS}, when $k \geq 4$ we have
    \begin{equation*}
        N 
        \leq MS_{k-2} \left( 2; n^{C' \binom{r}{s}^2 \log \binom{r}{s}} \right) 
        \leq T_{k-4}\left(2^{n^{C' \binom{r}{s}^2 \log \binom{r}{s}}}\right)
        = T_{k-3} \left( n^{C' \binom{r}{s}^2 \log \binom{r}{s}} \right).
    \end{equation*}
\qed


\section{A characterization of $A_k(n;r,s)$}\label{scharacterization}

\subsection{Generalized $S$-increasing sequences}

We will give a characterization of $A_k(n; r,s)$ in terms of what we call ``generalized $S$-increasing sequences" which simultaneously generalizes Loh's characterization for the 2-uniform problem (\cite[Lemma 2.1]{Loh15}) and Moshkovitz and Shapira's characterization for their problem in all uniformities $k$ (\cite[Theorem 5]{MS14}).
These generalized $S$-increasing sequences extend the (2-uniform) notion of $S$-increasing sequences, which were studied extensively by Gowers and Long in \cite{GL21}.

\begin{definition}\label{Sincreasing}
    Let $S < R$ and $n$ be positive integers.
    We define $M_k(n;R,S)$ and $m_k(n;R,S)$ by recursion on $k$ as follows.
    \begin{itemize}
        \item For $\vx = (x_1, x_2, \dots, x_R), \vy = (y_1, y_2, \dots, y_R) \in [n]^R$, we say that $\vx$ is \textit{$S$-less} than $\vy$, denoted $\vx <_S^{(2)} \vy$, if $x_\l < y_\l$ for at least $S$ coordinates $\l \in [R]$.
        A sequence $(\vx_i)$ in $[n]^r$ is \textit{$S$-increasing} if $\vx_i <_S^{(2)} \vx_j$ for all $i < j$.
        Let $M_2(n;R,S)$ be the set of all $S$-increasing sequences in $[n]^R$, and let $m_2(n;R,S)$ denote the length of the longest sequence in $M_2(n;R,S)$.
        \item For $k \geq 2$, let $\mathcal{P}^{(k-2)}([n]^R)$ denote the $(k-2)$-fold power set $\mathcal{P} ( \dots ( \mathcal{P}([n]^R)))$ of $[n]^R$, with $\mathcal{P}^{(0)}([n]^R) = [n]^R$, and recursively assume that we have defined a relation $<_S^{(k)}$ on $\mathcal{P}^{(k-2)}([n]^R)$.
        We define the relation $<_S^{(k+1)}$ on $\mathcal{P}^{(k-1)}([n]^R)$ as follows: 
        for two sets $\mathcal{F}, \mathcal{F}' \in \mathcal{P}^{(k-1)}([n]^R)$, we put $\mathcal{F} <_S^{(k+1)} \mathcal{F}'$
        if there exists a distinguished element $F' \in \mathcal{F}'$ such that $F <_S^{(k)} F'$ for all $F \in \mathcal{F}$.
        We will set $M_{k+1}(n;R,S)$ to be the set of all sequences $(\mathcal{F}_i)$ in $\mathcal{P}^{(k-1)}([n]^R)$ with the property that for all $i < j$, $\mathcal{F}_i <_S^{(k+1)} \mathcal{F}_j$. 
        Sequences in $M_{k+1}(n;R,S)$ will also be called \textit{$S$-increasing}, and we will denote by $m_{k+1}(n;R,S)$ the length of the longest sequence in $M_{k+1}(n;R,S)$.
    \end{itemize}
\end{definition}

\begin{example}
Here are explanations in words for $k=3$ and $k=4$ of what is meant in Definition~\ref{Sincreasing}:
    \begin{itemize}
        \item Elements of $M_3(n;R,S)$ are sequences $(F_1, F_2, \dots)$ of \textit{subsets} $F_i \subseteq [n]^R$ such that for all $i < j$, $F_i <_S^{(3)} F_j$, i.e., there exists $\vx' \in F_j$ with $\vx <_S^{(2)} \vx'$ for every $\vx \in F_i$;
        \item Elements of $M_4(n;R,S)$ are sequences $(\cF_1, \cF_2, \dots)$ of \textit{families of subsets} $\cF_i \subseteq 2^{[n]^R}$ with the property that for all $i < j$, $\cF_i <_S^{(4)} \cF_j$, i.e., there is some $F' \in \cF_j$ for which $F <_S^{(3)} F'$ holds for all $F \in \cF_i$.
    \end{itemize}
\end{example}

Our next theorem shows that $A_k(n;r,s)$ is always upper bounded by $m_k(n;\binom{r}{s},\binom{r-1}{s-1}) + 1$, and moreover that this upper bound is sharp for both the usual Ramsey regime $s=1$ (i.e., $MS_k(n;r)$) and the Erd\H{o}s-Szemer\'edi color-avoiding regime $s=r-1$.

\begin{theorem}\label{characterization}
    For all integers $n \geq k \geq 2$ and $r > s \geq 1$, $A_k(n; r,s) \leq m_k(n;\binom{r}{s}, \binom{r-1}{s-1}) + 1$.
    When $s \in \{1,r-1\}$, the lower bound $A_k(n;r,s) > m_k(n; \binom{r}{s}, \binom{r-1}{s-1})$ also holds.
\end{theorem}

\begin{proof}[Proof that $A_k(n; r, s) \leq m_k(n;r,s) + 1$]
    Let $\Delta: [N]^{(k)} \to [r]$ be a coloring in which every tight monotone path on $n$ edges uses at least $s+1$ colors.
    Define maps $P_{k-1}, P_{k-2}, \dots, P_1$ by reverse-recursion as follows:
    \begin{itemize}
        \item[(1)] Let $P_{k-1}: [N]^{(k-1)} \longrightarrow [n]^{\binom{r}{s}}$ for $v_2 < \dots < v_k$ be given by
        \begin{equation*}
            P_{k-1}(v_2, \dots, v_k) = \langle 1+\l_S(v_2, \dots, v_k) \rangle_{S \in [r]^{(s)}} \in [n]^{\binom{r}{s}},
        \end{equation*}
        where $\l_S(v_2, \dots, v_k)$ denotes the longest monotone path ending at $v_2 < \dots < v_k$ which only uses colors in the set $S \in [r]^{(s)}$;
        \item[(2)] For $2 \leq i \leq k-1$, define $P_{k-i}: [N]^{(k-i)} \longrightarrow \mathcal{P}^{(i-1)}([n]^{\binom{r}{s}})$ for $v_{i+1} < \dots < v_k$ by 
        \begin{equation*}
            P_{k-i}(v_{i+1}, \dots, v_k) = \{P_{k-i+1}(u, v_{i+1}, \dots, v_k) \in \mathcal{P}^{(i-2)}([n]^{\binom{r}{s}}): u < v_{i+1}\}.
        \end{equation*}
    \end{itemize}
    Notice that $P_1$ is a map from $[N]$ into $\mathcal{P}^{(k-2)}([n]^{\binom{r}{s}})$.
    We claim that the sequence $(P_1(v))_{v=1}^N$ is in  $M_k(n;{r\choose s},{r-1 \choose s-1})$ which would give (setting $N = A_k(n+1;r,s) - 1$)
    \begin{equation*}
        A_k(n;r,s) - 1 \leq m_k \left( n;\binom{r}{s},\binom{r-1}{s-1} \right).
    \end{equation*}
    To see this, suppose that $u < v$ are in $[N]$.
    We need to show $P_1(u) <_{\binom{r-1}{s-1}}^{(k)} P_1(v)$.
    First, note that for any $s_1 < \dots < s_{k-2} < u$, we have
        \begin{equation*}
            P_{k-1}(s_1, \dots, s_{k-2}, u) <_{\binom{r-1}{s-1}}^{(2)} P_{k-1}(s_2, \dots, s_{k-2}, u, v)
        \end{equation*}
    in $[n]^{\binom{r}{s}}$, as the longest paths ending at $(s_1, \dots, s_{k-2}, u)$ which only use colors in $S \in [r]^{(s)}$ with $\Delta(s_1, \dots, s_{k-2}, u, v) \in S$ can be extended one edge longer to end at $(s_2, \dots, s_{k-2}, u, v)$ while still only using colors in $S$; 
    as such, the longest paths ending at $(s_2, \dots, s_{k-2}, u, v)$ only using colors in these $\binom{r-1}{s-1}$ sets $S$ will be strictly longer than those ending at $(s_1, \dots, s_{k-2}, u)$.

    Now assume inductively that for $1 \leq i \leq k-1$, we have shown
    \begin{equation*}
        P_{k-i}(s_i, \dots, s_{k-2}, u)~<_{\binom{r-1}{s-1}}^{(i+1)}~P_{k-i}(s_{i+1}, \dots, s_{k-2}, u, v)
    \end{equation*}
    holds for any choice of $s_i < s_{i+1} < \dots < s_{k-2} < u$.
    Since 
    $$P_{k-i}(s_{i+1}, \dots, s_{k-2}, u, v) \in P_{k-(i+1)}(s_{i+2}, \dots, s_{k-2}, u,v)$$ 
    and every $F \in P_{k-(i+1)}(s_{i+1}, \dots, s_{k-2}, u)$ has the form $F = P_{k-i}(s_i, s_{i+1}, \dots, s_{k-2}, u)$ for some $s_i < s_{i+1}$, we get that 
    $$P_{k-(i+1)}(s_{i+1}, \dots, s_{k-2}, u) <_{\binom{r-1}{s-1}}^{(i+2)} P_{k-(i+1)}(s_{i+2}, \dots, s_{k-2}, u,v)$$
    with distinguished element $F' = P_{k-i}(s_{i+1}, \dots, s_{k-2}, u, v)$ immediately from our inductive hypothesis. To finish, setting $i = k-1$ gives $P_1(u) <_{\binom{r-1}{s-1}}^{(2)} P_1(v)$, as desired.
\end{proof}

\textit{Proof that $A_k(n; r,s) > m_k(n; \binom{r}{s}, \binom{r-1}{s-1})$ for $s \in \{1,r-1\}$.}
    Let $m_k(n; \binom{r}{s}, \binom{r-1}{s-1}) = N$ and let $(\mathcal{F}_i)_{i=1}^{N}$ be a maximum-length sequence in $M_k(n; \binom{r}{s}, \binom{r-1}{s-1})$.
    We define an $r$-coloring $\Delta: [N]^{(k)} \to [r]$ in which no tight monotone path on $n$ edges uses at most $s$ colors.\\
    Let $i_1 < i_2 < \dots < i_k$ be elements of $[N]$, and consider $\mathcal{F}_{i_1}, \dots \mathcal{F}_{i_k}$.
    Since the sequence of $\mathcal{F}_i$'s is $\binom{r-1}{s-1}$-increasing, we know that for each $j \in \{2, \dots, k\}$, there exists $F_{i_{j-1}, i_j} \in \mathcal{F}_{i_j}$ such that $F <_{\binom{r-1}{s-1}}^{(k)} F_{i_{j-1}, i_j}$ for all $F \in \mathcal{F}_{i_{j-1}}$.
    In particular, we have
    \begin{equation}\label{relation}
        F_{i_1, i_2} <_{\binom{r-1}{s-1}}^{(k)} F_{i_2, i_3} <_{\binom{r-1}{s-1}}^{(k)} \dots <_{\binom{r-1}{s-1}}^{(k)} F_{i_{k-1}, i_k}.
    \end{equation}
    Now repeat this $k-2$ more times | each time, we find distinguished elements in all but the first set which satisfy \eqref{relation}.
    At the end we are left with $\vx_{i_1, \dots, i_{k-1}} <_{\binom{r-1}{s-1}}^{(2)} \vx_{i_2, \dots, i_k}$ in $[n]^{\binom{r}{s}}$.

    \begin{itemize}
        \item If $s=1$: 
        Then $\binom{r-1}{s-1} = 1$, and so $\vx_{i_1, \dots, i_{k-1}} <_1^{(2)} \vx_{i_2, \dots, i_k}$, i.e., there is some coordinate in which $\vx_{i_1, \dots, i_{k-1}}$ is strictly less than $\vx_{i_2, \dots, i_k}$, and we can define our coloring $\Delta$ by setting $\Delta(i_1, \dots, i_k)$ to be the least coordinate where this occurs.
        Then any monotone tight path on the $t$ edges on vertices $i_1 < \dots < i_{t+k-1}$ using at most one color $c \in [\binom{r}{1}] = [r]$ (i.e., which is monochromatic) would yield a sequence
        \begin{equation*}
            \vx_{i_1, \dots, i_{k-1}} <_1^{(2)} \vx_{i_2, \dots, i_k} <_1^{(2)} \dots <_1^{(2)} \vx_{i_{t+1}, \dots, i_{t+k-1}}
        \end{equation*}
        which strictly increases in the coordinate $c$, whence $t \leq n-1$.

        \item If $s=r-1$: Then $\binom{r-1}{s-1} = r-1$ and $\vx_{i_1, \dots, i_{k-1}} <_{r-1}^{(2)} \vx_{i_2, \dots, i_k}$, so in all but one coordinate the entry in $\vx_{i_1, \dots, i_{k-1}}$ is strictly less than the one in $\vx_{i_2, \dots, i_k}$.
        Therefore, we can define our coloring $\Delta$ by setting $\Delta(i_1, \dots, i_k)$ to be the (lexicographically) least $(r-1)$-set of coordinates on which $\vx_{i_1, \dots, i_{k-1}} <_{r-1}^{(2)} \vx_{i_2, \dots, i_k}$ holds.
        Now any monotone tight path on the $t$ edges on vertices $i_1 < \dots < i_{t+k-1}$ using at most $r-1$ colors yields a sequence
        \begin{equation*}
            \vx_{i_1, \dots, i_{k-1}} <_{r-1}^{(2)} \vx_{i_2, \dots, i_k} <_{r-1}^{(2)} \dots <_{r-1}^{(2)} \vx_{i_{t+1}, \dots, i_{t+k-1}}
        \end{equation*}
        where there are at most $r-1$ choices for the set of $r-1$ coordinates where consecutive elements strictly increase.
        Since this means that at most $r-1$ of the $\binom{r}{r-1} = r$ coordinates can be left out by one of these sets, there is some coordinate which is in all of them, and the $\vx$'s must strictly increase on this coordinate.
        Hence, we again have $t \leq n-1$. \qed
    \end{itemize}

\section{Concluding remarks}\label{conclusion}

There are several problems related to this paper which remain open.
The most immediate is, of course, to give good lower bounds for $A_k(n;r,s)$. 
To this end, we pose the following question.

\begin{question}
    For all $n$ sufficiently large, $k \geq 4$, and $r > s > r/2$, is it true that
    \begin{equation*}
        A_k(n;r,s) \geq T_{k-3}(n^c)
    \end{equation*}
    for some $c = c(r,s) > 0$?
\end{question}

In other words, is there a construction showing that the tower height from our Theorem \ref{maink} cannot be improved? 
We know that this is the case in uniformities $k=2$ and $k=3$, but we do not know what the truth is in uniformity $k \geq 4$.
By Theorem~\ref{characterization}, for $s=r-1$ it is equivalent to study the function $m_k(n; \binom{r}{s}, \binom{r-1}{s-1}) = m_k(n; r, r-1)$ introduced in Section~\ref{scharacterization}.
It is worth pointing out while the usual approach to obtain tower-type lower bounds for hypergraph Ramsey numbers, known as the \textit{stepping-up lemma} of Erd\H{o}s and Hajnal \cite{EHR65}, adapts very well for the problem of lower bounding the function $MS_{k}(n;r)$ (the lower bound in Theorem \ref{MS} above; see \cite{FPSS12, MS14}), but fails for interesting reasons for this problem.

The corresponding color-avoiding problem for cliques is also open. Stepping-up lemmas for various combinations of $r$, $s$, and $k$ were proved in~\cite{DGHY23}, but the quantitative aspects of the problem with $r=3$, $s=2$ are still poorly understood even in uniformity three.
In this regime, a first-moment method argument yields an exponential-type lower bound, while standard arguments give a double-exponential upper bound, so the tower height remains undetermined.
Based on how our result distinguishes $A_3(n;3,2)$ from $MS_3(n;2)$, we would like to make the following conjecture about the behavior of the three-uniform color-avoiding Ramsey function for $r=3$ and $s=2$:

\begin{conjecture}
    Let $R_{3}(n;3,2)$ be the least integer $N$ with the property that in every three-coloring of $E(K_N^{(3)})$, there exists a copy of $K_n^{(3)}$ which uses at most two colors.
    Then
    \begin{equation*}
        R_{3}(n;3,2) \leq 2^{n^{c}}
    \end{equation*}
    for some $c > 0$.
\end{conjecture}

Another potentially fruitful direction would be to explore the case when $s \leq r/2$ (for both the $A_{3}(n;3,2)$ and $R_{3}(n;3,2)$). In this regime, Theorem \ref{kmaj} no longer applies, so it is an interesting problem to determine whether or not $A_{3}(n;r,s)$ still exhibits a drastically different behavior than $MS_{3}(n;\lceil r/s \rceil)$, for example. 

Last but not least, we would also like to mention recent work of Suk and Zeng \cite{suk2024monotone} which establishes a connection between Ramsey numbers of 3-uniform monotone paths and the Erd\H{o}s-Schur problem.





\section*{Acknowledgments} 
The authors would like to thank the anonymous reviewers for several helpful suggestions.

\bibliographystyle{amsplain}


\begin{dajauthors}
\begin{authorinfo}[eion]
  Eion Mulrenin\\
  Department of Mathematics\\
  Emory University\\
  Atlanta, USA\\
  eion\imagedot{}mulrenin\imageat{}emory\imagedot{}edu \\
\end{authorinfo}
\begin{authorinfo}[cosmin]
  Cosmin Pohoata\\
  Department of Mathematics\\
  Emory University\\
  Atlanta, USA\\
  cosmin\imagedot{}pohoata\imageat{}emory\imagedot{}edu \\
  \url{https://pohoatza.wordpress.com/about/}
\end{authorinfo}
\begin{authorinfo}[dima]
  Dmitrii Zakharov\\
  Department of Mathematics\\
  Massachusetts Institute of Technology\\
  Cambridge, USA\\
  zakhdm\imageat{}mit\imagedot{}edu\\
  \url{https://sites.google.com/view/dmitriizakharov/}
\end{authorinfo}
\end{dajauthors}

\end{document}